\documentclass[3p,10pt,a4paper,twoside,fleqn,sort&compress]{filomat}
\usepackage{amssymb,amsmath,latexsym,color}
\usepackage[varg]{pxfonts}
\usepackage{dsfont}

\newtheorem{theorem}{Theorem}[section]

\newtheorem{corollary}[theorem]{Corollary}

\newtheorem{definition}[theorem]{Definition}

\newtheorem{lemma}[theorem]{Lemma}

\newtheorem{proposition}[theorem]{Proposition}

\newcommand{\FilVol}{xx}
\newcommand{\FilYear}{yyyy}
\newcommand{\FilPages}{zzz--zzz}
\newcommand{\FilDOI}{(will be added later)}
\newcommand{\FilReceived}{Received: xx-xx-xxxx; Accepted: xx-xx-xxxx}
\begin{document}

	\title{A convergence theorem for $ap-$Henstock-Kurzweil integral
		and its relation to topology	}
	\author[affil1]{Hemanta Kalita$^\ast$}
	\ead{hemanta30kalita@gmail.com; hk\_gu@gauhati.ac.in} 
	\author[affil2]{Bipan Hazarika}
	\ead{bh\_rgu@yahoo.co.in; bh\_gu@gauhati.ac.in}
	

	\address[affil1]{Department of Mathematics, Assam DonBosco University, Sonapur, Guwahati 782402, Assam, India}
	\address[affil2]{Department of Mathematics, Gauhati University, Gauhati 781014, Assam, India}

	\newcommand{\AuthorNames}{Hemanta Kalita, Bipan Hazarika}
	
	\newcommand{\FilMSC}{26A39, 46B03, 46B20, 46B25}
	\newcommand{\FilKeywords}{$ap-$Henstock-Kurzweil integrable function;  locally convex topology; topology in the primitive class}
	\newcommand{\FilCommunicated}{(Binod Chandra Tripathy)}
	\newcommand{\FilCorrespondingauthor}{$^\ast$Hemanta Kalita}
	\newcommand{\FilSupport}{(This article is dedicated to the memory of Prof. Henri Poinca\'re).}
	

 





\begin{abstract} In this {\color{red}{paper}} we discuss about the  $ap-$Henstock-Kurzweil integrable functions on a topological vector spaces. Basic results of $ap-$Henstock-Kurzweil integrable functions are discussed here. We discuss the  equivalence of  the $ap-$Henstock-Kurzweil integral on a topological vector spaces and  the vector valued $ap-$Henstock-Kurzweil integral. Finally, several convergence theorems are studied. 
\end{abstract}
\maketitle{}
\makeatletter
\renewcommand\@makefnmark%
{\mbox{\textsuperscript{\normalfont\@thefnmark}}}
\makeatother



\section{Introduction}
A further approach to the problem of the primitives was introduced in
1957 by J. Kurzweil  and in 1963 by R. Henstock, independently.
They defined a generalized version of the Riemann integral that is known as
the Henstock-Kurzweil integral, also abbreviated as the HK-integral. The
advantage of the HK-integral is that, it is very similar in construction and
in simplicity to the Riemann integral and it has the power of the Lebesgue
integral. Moreover, in the real line, the HK-integral solves the problem of
the primitives. The definition of the HK-integral is constructive, as in the
Riemann integral, and the value of the HK-integral is defined as the limit
of Riemann sums over suitable partitions of the domain of integration. The
main difference between the two definitions is that, in the HK-integral, a
positive function, called gauge, is used, instead of the constant utilized in
the Riemann integral to measure the fineness of a partition (one can see \cite{B,HS,HI,J,Base,KW,RA}). This gives a
better approximation of the integral near singular points of the function. For integration of approximate derivative the situation turned out to be more complicated. Most of researchers effort in this field was exerted into finding relations between approximate Perron-type integrals and the Denjoy-Khintchine integral and its approximately continuous generalizations. The approximately continuous Perron integral (AP-integral) was introduced by Burkill  \cite{Burkill}.  Park et al. \cite{Jae1} studied  the convergence theorem for the AP-integral based on the condition UAP and pointwise boundedness. Park et al. \cite{Jae} defined the AP-Denjoy integral and show that the AP-Denjoy
integral is equivalent to the AP-Henstock-Kurzweil  integral and the integrals are equal.  Skvortsov and Sworowski  \cite{Val1} brought to attention on the  known results which are stronger than those contained in the  work of \cite{Jae}. They show that some of them can be formulated
in terms of a derivation basis defined by a local system of which the approximate basis
is known to be a particular case. They also consider the relation between the $\sigma-$finiteness
of variational measure generated by a function and the classical notion of the generalized
bounded variation. For a wide class of bases, Riemann-type integral is equivalent to the appropriately defined Perron-type integral (see \cite{O}).  Skvortsov et al. \cite{Val}  say only that Burkill’s $ap-$integral is covered by $ap-$Henstock-Kurzweil integral.  Shin et al.  \cite{Kwan}  introduced the concept of approximately negligible variation and give a necessary and sufficient condition that a function $\mathcal{F}$ be an indefinite
integral of an $ap-$Henstock integrable function $f$ on $[a, b].$ They characterized absolutely $ap-$Henstock integrable functions by using the concept of bounded variation. Yoon \cite{JU}, studied  about vector valued $ap-$Henstock-Kurzweil integrals. They discussed  some of its properties, and characterize $ap-$Henstock integral of vector valued functions by the notion of equiintegrability. It is known that for Banach valued function Henstock-Lemma fail for Henstock-Kurzweil integral but Henstock-Lemma holds for a locally convex spaces (see \cite{VM}). Yoon \cite{JU},  has not discussed about Henstock-Lemma for their integrals. This is an open area till now. We motivate from the article \cite{VM}, that vector valued concept of the Henstock-Kurzweil integral is not at all sufficient for overall studies of the Henstock-Kurzweil integral. We introduce  the concept of $ap-$Henstock-Kurzweil integrals on topological vector spaces and investigates several convergence theorems in this settings.
\section{Preliminaries}
Let $X$ be a Hausdorff topological space. We say that $X$ is a topological
vector space (in short TVS) if $X$ is a real vector space and the operations, vector addition and
scalar multiplication, are continuous.
\begin{definition}
Let $X$ be a non empty set. A family $\mathfrak{F}=\{A_\nu:~\nu \in \mathbb{N}\}$ of subsets of $X$ is a filter in $X$ if the following are satisfied:
\begin{enumerate}
\item For every $\nu \in \mathbb{N},~A_\nu \neq \emptyset.$
\item For $A,B \in \mathfrak{F}$ then $A \cap B \in \mathfrak{F}.$
\item If $A \in \mathfrak{F},~B \subseteq X$ and $A \subseteq B$ then $B \in \mathfrak{F}.$
\end{enumerate}
\end{definition}
   The filter $\mathfrak{F}$ converges to $ x \in X$ if for every $\theta-$nbd $U$ ($\theta$ is zero vector of $X$)  there exists $A \in \mathfrak{F}$ such that $A - x \subseteq U.$  We say $\mathfrak{F}$ is Cauchy if for every $\theta-$nbd $U$ there exists $A \in \mathfrak{F}$ such that $A - A \subseteq U.$
   \begin{definition}
   Given a measurable set $E \subset [a,b],$ a set valued function $\Delta: E \to 2^{[a,b]}$ is an $ap~\theta-$nbd function (ANF) on $E$ if for every $x \in E,$ there exists an $ap~\theta-$nbd $U_x \subset [a,b]$ of $x$ such that $\Delta(x)=U_x.$
   \end{definition}
   \begin{definition}
   Let $f:[a,b] \to X,~F:[a,b] \to X$ and let $E \subset [a,b]$ be a measurable. $F$ is said to satisfy the approximate strong Lusin conditions on $E( F \in ASL(E))$ if for every $Z \subset E$ of measure zero and for every $\varepsilon>0$ there exists an ANF $\Delta$ on $E$ such that $$S(|F|,P)- \mathcal{A}) \in U$$ for a $\theta-$nbd $U.$
   \end{definition}
    Let $\mathbb{I}$ denote all non degenerated closed intervals of $[a,b]$ and $\lambda$ be the Lebesgue measure on $[a,b].$ We denote an interval function $F:\mathbb{I} \to \mathbb{R}$ with the end point $F(t)=F([a,t]),~t \in [a,b].$ That is, $F([e,f])= F(f)-F(e),~[e,f] \in \mathbb{I}.$ Throughout the {\color{red}{paper}} measurable functions are mean by $\lambda-$measurable.\\ Recalling when $X$ is a Banach space the $ap-$Henstock-Kurzweil integral is as follows
   \begin{definition}\label{banach}
  \cite[Definition 2.1]{JU}  A function $f : [a, b] \to  X$ is $ap-$Henstock integrable on $[a, b]$ if there exists a vector $\mathcal{A} \in X$  with the following property:\\
for each $\varepsilon >0$ there exists a choice $S$  on $[a, b]$ such that $||S(f, P)-\mathcal{A}||< \varepsilon$
whenever $P$ is a tagged partition of $[a, b]$ that is subordinate to $S.$ The
vector $\mathcal{A}$  is called the $ap-$Henstock integral of $f$ on $[a, b]$ and is denoted by $(ap)\int_{a}^{b}f.$
   \end{definition}
   \section{ Basic properties of $ap-$Henstock-Kurzweil integral on a Topological Vector Spaces}
   An approximate $\theta-$nbd of $x \in [a,b]$ is a measurable set $S_x \subset [a,b]$ containing $x$ as a point of density. Let $E \subset [a,b].$ For every $x \in E \subset [a,b],$ choose an approximate $\theta-$nbd $S_x \subset [a,b] $ of $x.$ Then $S=\{S_x:~x \in E\}$ is a choice on $E.$ We assume each point of $S_x$ is a point of density of $S_x.$\\ A tagged interval $([u,v],x)$ is said to be fine to the choice $S=\{S_x\}$ if $u,v \in S_x$ and $x \in [u,v].$ A tagged sub partition $P=\{([u_i,v_i],t_i):~1 \leq i \leq n\}$ of $[a,b]$ is a finite collection of non overlapping tagged interval in $[a,b]$ such that $t_i \in [u_i, v_i]$ for $i=1,2,\dots,n,$ then we say $P$ is $S-$fine. If $P$ is $S-$fine and $t_i \in E$ for each $1 \leq i \leq n,$ then $P$ is $(S,E)-$fine. If $P$ is $S-$fine and $[a,b]= \bigcup_{i=1}^{n}[u_i,v_i],$ then we say $P$ is $S-$fine tagged partitions of $[a,b].$
   \begin{definition}
   \begin{enumerate}
   \item[(1)] $f:[a,b] \to \mathbb{R}$ is approximately continuous at $c \in [a,b]$ if there exists a measurable $\theta-$nbd $U \subset [a,b]$ with density $1$ at $c$ such that $f(x)-f(c) \in U$ whenever $|x-c|< \delta.$
   \item[(2)] We say $f$ is approximately differentiable at $c$ if there exists a real number $A$ and a measurable $\theta-$nbd $U \subset [a,b]$ such that the density of $U$ at $c$ is $1$ and $\frac{f(x)-f(c)}{x-c}-A \in U.$
   \end{enumerate}
   \end{definition}
   For a tagged partition $P=\{([u_i,v_i],t_i):~ 1 \leq i \leq n \}$ of $[a,b]$ we define the Riemann sum as $$S(f,P) = \sum_{i=1}^{n}f(t_i)(v_i - u_i)~\mbox{~if~it exists}.$$
   \begin{definition}\label{tvs}
   A function $f:[a,b] \to X$ is $ap-$Henstock-Kurzweil integrable on $[a,b]$ if there exists an $\mathcal{A} \in X$ such that for any $\theta-$nbd  $U$ of $[a,b]$  there exists a  gauge $\delta$ on $[a,b]$ whenever $$P=\{([x_{i-1},x_i],t_i):~ 1 \leq i \leq n \}$$ is $S-$fine of $[a,b],$ we have $$S(f,P)- \mathcal{A} \in U.$$  We call $\mathcal{A}$ is the $ap-$Henstock-Kurzweil integral of $f$ on $[a,b].$ Here $\mathcal{A}= (ap)\int_{a}^{b}f .$
   \end{definition}
   Let us consider $AP([a,b], X)$ be the set of all $ap-$Henstock-Kurzweil integrable $X-$valued functions on $[a,b].$ The function $f$ is $ap-$Henstock-Kurzweil integrable on a measurable set $E \subseteq [a,b]$ if $f \chi_E$ is $ap-$Henstock-Kurzweil integrable on $[a,b],$ where $\chi_E$ is the characteristic functions on $E.$ In this settings the Henstock-Kurzweil integrable functions are certainly the $ap-$Henstock-Kurzweil integrable.
   \begin{definition}
   If $f \in AP([a,b], X)$ and $u \in [a,b],$ then $F(u)= \int_{a}^{u}f$ is called $ap-$primitive of the $ap-$Henstock-Kurzweil integral $f.$
   \end{definition}
   If $P=\{([x_{i-1},x_i],t_i)\}_{i=1}^{n}$ is any $S-$fine tagged partition on $[a,b]$ then we denote $$F(P)= \sum_{i=1}^{n}F([x_{i-1},x_i]).$$
   \begin{proposition}
   Every given function $f:[a,b] \to X$ have at most one $ap-$Henstock-Kurzweil integral on $[a,b].$
   \end{proposition}
   \begin{proof}
   Suppose $f$ is $ap-$Henstock-Kurzweil integrable on $[a,b].$ If possible, let us consider $\mathcal{A}_1$ and $\mathcal{A}_2$ are the $ap-$Henstock-Kurzweil integral of $f$ with $\mathcal{A}_1 \neq \mathcal{A}_2.$ Let $U_1$ and $U_2$ be disjoint $\theta-$nbds of $\mathcal{A}_1$ and $\mathcal{A}_2,$ respectively. The fact from the $\theta-$nbd, $U_1 - \mathcal{A}_1$ and $U_2 - \mathcal{A}_2$ are $\theta-$nbd. Say $W_1= U_1 - \mathcal{A}_1$ and $W_2= U_2 - \mathcal{A}_2,$ then for every $S-$fine tagged partition $P_1$ and $P_2$ of $[a,b]$ with $P_1$ is $\delta_1-$fine and $P_2$ is $\delta_2-$ fine tagged on $[a,b]$  we have 
   \begin{eqnarray}
   &&S(f,P_1)- \mathcal{A}_1 \in W_1 \\&& S(f,P_2)- \mathcal{A}_2 \in W_2.
   \end{eqnarray}
   Let $\delta(x)= \min( \delta_1(x), \delta_2(x))$ for all $x \in [a,b]$ and $P$ be a $S-$fine tagged partitions of $[a,b].$ Clearly $P$ is $\delta_1,~\delta_2-$fine. Hence, $S(f,P) \in W_1 + \mathcal{A}_1 = U_1$ and $S(f,P) \in W_2 + \mathcal{A}_2 = U_2.$ This is a contradiction. So, $\mathcal{A}_1 = \mathcal{A}_2.$
   \end{proof}
   \begin{proposition}
   If $X$ be a topological vector space. If $\alpha$ is a real number and $f,g \in AP([a,b], X)$ then $\alpha f,~f+g \in AP([a,b], X)$ with $$(ap)\int_{a}^{b}\alpha f = \alpha (ap)\int_{a}^{b} f $$ and $$(ap)\int_{a}^{b}(f+g)= (ap)\int_{a}^{b}f +(ap)\int_{a}^{b} g.$$
   \end{proposition}
   \begin{proof}
   Let us assume $(ap)\int_{a}^{b} f = \mathcal{A}.$ If $\alpha=0,$ then $(ap)\int_{a}^{b}\alpha f = \alpha (ap)\int_{a}^{b} f .$ Assume $\alpha \neq 0,$ let $U$ be $\theta-$nbd then there exists a $S-$fine partition $P$ of $[a,b]$ such that $$S(f,P)-A \in \frac{U}{\alpha}.$$ Thus,
   \begin{align*}
   S(\alpha f, P) -\alpha \mathcal{A} &= \alpha(S(f,P))- \alpha f\\&= \alpha(S(f,P)- \mathcal{A})\\
   &\in U.
   \end{align*}
   Hence $\alpha f \in AP([a,b],X)$ and $(ap)\int_{a}^{b}\alpha f = \alpha (ap)\int_{a}^{b} f .$\\ For the second part, let $(ap)\int_{a}^{b}f = \mathcal{A}_1$ and $(ap)\int_{a}^{b} g = \mathcal{A}_2.$ If  $U$ be $\theta-$nbd then there exists a $\theta-$nbd $V$ (say) such that $V+V \leq U.$ Consequently, $S(f,P_1) - \mathcal{A}_1 \in V$ for a $S-$fine partition $P_1$ on $[a,b].$ With the similar fashion for a $S-$fine tagged partition $P_2$ on $[a,b],$ we have $S(f,P_2) - \mathcal{A}_2 \in V.$ If $P = \min(P_1, ~P_2),$ then we have $$S(f+g,P)= S(f,P)+S(g,P).$$ Thus,
   \begin{align*}
   S(f+g,P) -(\mathcal{A}_1+ \mathcal{A}_2)&= S(f,P)- \mathcal{A}_1 + S(g,P)- \mathcal{A}_2\\& \in V+V \subseteq U.
   \end{align*}
   Hence $f+g \in AP([a,b],X)$ and $(ap)\int_{a}^{b}f+g = (ap)\int_{a}^{b}f + (ap)\int_{a}^{b}g.$
   \end{proof}
   \begin{proposition}
   Let $X$ be a topological vector space. If $f \in AP([a,b],X)$ and $f \in AP([b,c],X)$ then $f \in AP([a,c],X)$ and $$(ap)\int_{a}^{c}f= (ap)\int_{a}^{b}f + (ap)\int_{b}^{c}f.$$
   \end{proposition}
   \begin{proposition}
   (Cauchy's criterion) Let $X$ be a complete topological vector space. Then $f \in AP([a,b],X)$ if and only if for every $\theta-$nbd $U$ there exists a $S-$fine gauge $\delta$ on $[a,b]$ such that $$S(f,P_1)-S(f,P_2) \in U$$ for each pair $S-$fine partitions $P_1$ and $P_2$ of $[a,b].$
   \end{proposition}
   \begin{proof}
   Let us assume $(ap)\int_{a}^{b}f= \mathcal{A}.$ If $U$ is a $\theta-$nbd then there exists a $\theta-$nbd $V$ such that $V-V \subseteq U.$ From the definition of the $ap-$Henstock-Kurzweil integral $S(f,P) - \mathcal{A} \in V$ for a $S-$fine partition $P$ of $[a,b].$ Now for $P_1,~P_2$ as $S-$fine partitions of $[a,b],$ we have
   \begin{align*}
   S(f,P_1) -S(f,P_2) &= (S(f,P)-\mathcal{A}) ) - (S(f,P_2) - \mathcal{A}_2)\\& \in V- V \subseteq U.
   \end{align*}
   Let $A_\delta = \{S(f,P): P\mbox{ ~is~S-fine~tagged~partition~of ~} [a,b]\}.$ Also, assume $$\mathbb{A}=\{A_\delta:~\delta-{\mbox{is~S-fine~tagged~on}} [a,b]\}.$$ Then clearly $\mathbb{A}$ is filter base in $X.$ From the completeness of $X$ we get $ \mathbb{A} \to \mathcal{A}$ for some $\mathcal{A} \in X.$ If $\mathcal{A}$ is $ap-$Henstock-Kurzweil integrable then our claim will over. Since $\mathbb{A} \to \mathcal{A}$ then $A_\delta - \mathcal{A} \subseteq U.$ Thus if $P$ is $S-$fine partition on $[a,b]$ then we have $S(f,P) -\mathcal{A} \in U.$ So, $f \in AP([a,b],X).$
   \end{proof}
   \begin{theorem}\label{final}
   Let $X$ be a complete topological vector space. A function $f:[a,b] \to X$ is in $AP([a,b], X)$ if and only if  the  following conditions are assure: \\For each $\theta-$nbd $U$ there exists a $S-$fine gauge $\delta$ on $[a,b]$ such that if $P=\{([x_{i-1}, x_i], t_i)\}_{i=1}^{n}$ is a $S-$fine tagged partition of $[a,b],$ also there exist open sets $U_i:~i=1,2,..,n$ with $\sum_{i=1}^{n}U_i \subseteq U$ and a function $F:[a,b] \to X$ such that 
   \begin{eqnarray}\label{hk1}
   F(x_i)-F(x_{i-1})-(x_i- x_{i-1})f(t_i) \in U_i~for~all~i.
   \end{eqnarray}
   \end{theorem}
   \begin{proof}
   We assume $f \in AP([a,b],X)$ and $F(x)= (ap) \int_{a}^{x}.$ If $U$ is $\theta-$nbd then there exists a $S-$fine gauge $\delta,$ a $S-$fine tagged partition $P=\{(x_{i-1},x_i],t_i)\}_{i=1}^{n}$ of $[a,b]$ such that $F(b)-S(f,P) \in U.$ Now, $$F(b)-S(f,P)= \sum_{i=1}^{n}(F(x_i)-F(x_{i-1})-(x_i - x_{i-1})f(t_i)).$$ Now from the given fact $\sum\limits_{i=1}^{n}U_i \subseteq U,$ clearly we get $$F(b)-S(f,P) \in U_i~\forall~i.$$ Conversely,  $f:[a,b] \to X$ assure  the equation (\ref{hk1}) for each $\theta-$nbd $U$ there exists a $S-$fine gauge $\delta$ on $[a,b]$ such that if $P=\{([x_{i-1}, x_i], t_i)\}_{i=1}^{n}$ is a $S-$fine tagged partition of $[a,b],$ also there exist open sets $U_i:~i=1,2,..,n$ with $\sum_{i=1}^{n}U_i \subseteq U.$ For,
   \begin{align*}
   F(b)-F(a)-S(f,P) &= \sum_{i=1}^{n} (F(x_i)-F(x_{i-1}-(x_i - x_{i-1})f(t_i))-F(a)\\&= \sum_{i=1}^{n} (F(x_i)-F(x_{i-1}-(x_i - x_{i-1})f(t_i))\\& \in \sum_{i=1}^{n}U_i \subseteq U.
   \end{align*}
   Hence $f \in AP([a,b],X).$
   \end{proof}
  \begin{theorem}
  (Saks-Henstock Lemma) Let $f:[a,b] \to X$ be $ap-$Henstock-Kurzweil integrable on $[a,b]$ and $F(x)= (ap)\int_{a}^{x} f$ for each $x \in [a,b],~\varepsilon>0.$ Suppose $S$ is a choice on $[a,b]$ such that for a $\theta-$nbd $U$ gives $S(f,P)-F(P) \in U.$ If $P_0=\{(J_j,t_j)\}_{j=1}^{s}$ is any $S-$fine tagged sub-partition of $[a,b]$ then $$S(f,P_0)-F(P_0) \in U$$ Moreover $$\sum_{j=1}^{s}(S(f,P_0)-(ap)\int_{J_j}f) \in U.$$ 
  \end{theorem}
  \begin{proof}
  The proof of the result follows from the Theorem \ref{final}.
  \end{proof}
   Now we discuss the topological vector space valued $ap-$Henstock-Kurzweil integrals are equivalent as the Banach valued $ap-$Henstock-Kurzweil integrals that discussed by Ju Han Yoon  \cite{JU}.
   \begin{theorem}
   Let $(X,||.||)$ be a Banach space. Then the $ap-$Henstock-Kurzweil integral are equivalents to the $ap-$Henstock-Kurzweil integrals for a topological vector spaces. That is, the Definition \ref{banach} and the Definition \ref{tvs} are equivalent.
   \end{theorem}
   \begin{proof}
   The Definition \ref{banach} implies the Definition \ref{tvs}: Let $U$ be $\theta-$nbd, then there exists $\varepsilon>0$ such that $B_\varepsilon \subseteq U,$ where $$B_\varepsilon = \{ x \in X: ||x||< \varepsilon \}.$$ Suppose $f$ satisfies the Definition \ref{banach}, then there exists a vector $\mathcal{A} \in X$ with the following property:: for each $\varepsilon>0$ there exists a choice $S$ on $[a,b]$ such that $$||S(f,P)-\mathcal{A}||< \varepsilon$$ whenever $P$ is $S-$fine tagged partition of $[a,b].$ Thus $$S(f,P) - \mathcal{A} \in B_\varepsilon \subseteq U.$$ Therefore $f$ satisfies the Definition \ref{tvs}.\\ Conversely, assume $f$ satisfies the Definition \ref{tvs}. Let $\varepsilon>0,$ then for a $S-$fine partition $P$ on $[a,b]$ we have $$S(f,P) -\mathcal{A} \in B_\varepsilon.$$ Thus $||S(f,P)- \mathcal{A}||< \varepsilon$ whenever $P$ is $S-$fine partition on $[a,b].$ This completes the  proof.
   \end{proof}
   \section{Convergence theorem on $ap-$Henstock-Kurzweil integrals and its relation to topology}
   In the literature the Denjoy convergence theorem generalizes the Vitali Convergence Theorem. The Perron convergence theorem  generalizes the Lebesgue Dominated Convergence Theorem. The $ap-$Henstock-Kurzweil convergence theorem generalized Dominated convergence theorem,  also the convergence theorem for the $ap-$Henstock-Kurzweil  integral based on the condition $UAP$ and pointwise boundedness. Here we study the  Dominated Convergence Theorem for the $ap$Henstock-Kurzweil integral and the convergence theorem for the $ap$Henstock-Kurzweil integral based on the condition uniformly $ap$Henstock-Kurzweil integrals  and the pointwise boundedness.
   \begin{definition}
Let $f: [a,b] \to X$ be measurable function. Let $\{f_k\}$ be a sequence of  integrable function defined on $[a,b].$ The sequence $\{f_k\}$ is said to be $ap-$Henstock-Kurzweil equi-integrable on $[a,b]$ if $\{f_k\}$ is $ap-$Henstock-Kurzweil integrable on $[a,b]$ if for each $\varepsilon>0$ there exists a choice $S$ such that $$S(f_k,P)-(ap)\int_{a}^{b} f_k d\mu \in U$$ hold for each $S-$fine partition $P=\{([x_{i-1}, x_i], t_i)\}_{i=1}^{n}$ of $[a,b],$  a $\theta-$nbd $U$ and $n \in \mathbb{N}.$
\end{definition}
 It is observe that if 
 $(f_n)$ be a pointwise bounded sequence of function $f_n:[a,b] \to X$ and let $E$ be subset of $[a,b]$ such that $\lambda(S\setminus E)=0.$ Then the sequence $(f_n)$ is $ap-$Henstock-Kurzweil integrable if and only if $f_n. \chi_E$ is $ap-$Henstock Kurzweil integrable.
   \begin{theorem}
Let $\{f_k\}_k$ be a non-decreasing sequence of $ap-$Henstock-Kurzweil integrable functions on  $[a,b]$ and let $f= \lim\limits_{k}f_k.$ If $\lim\limits_{k \to \infty}(ap)\int_{a}^{b}f_k  < \infty$ then $f$ is $ap-$Henstock-Kurzweil integrable on $[a,b]$ and $(ap)\int_{a}^{b}f = \lim\limits_{k \to \infty}(ap)\int_{a}^{b}f_k .$
\end{theorem}
\begin{proof}
From the definition of $ap-$Henstock-Kurzweil equi-integrability of $\{f_k\},$ for each $\varepsilon>0$ there exists a choice $S$ and a $\theta-$nbd $U$  such that \begin{eqnarray*}
S(f_k,P)-(ap)\int_{a}^{b} f_k \in U
\end{eqnarray*}
for each $S-$fine partition $P=\{([x_{i-1}, x_i], t_i)\}_{i=1}^{n}$ of $[a,b]$ and $n \in \mathbb{N}.$ Let $P$ be fixed. Since $\lim\limits_{n \to \infty}f_k(x)=f(x)$ then there is $m_0 \in \mathbb{N}$ such that 
\begin{eqnarray*}
S(f_k,P)-S(f_m,P) \in U~\forall~k, m >m_0.
\end{eqnarray*}
This implies $(ap)\int_{a}^{b}f_k  -(ap)\int_{a}^{b}f_m \in U$, therefore $(ap)\int_{a}^{b}f_k  $ of elements of $[a,b]$ is Cauchy and $\lim\limits_{k \to \infty}(ap)\int_{a}^{b}f_k  = \mathcal{A} \in X$ exists. This implies $$(ap)\int_{a}^{b}f_k  -\mathcal{A} \in U$$ for $m_1 \in \mathbb{N}$ with $k>m_1.$ Let any $S-$fine partition $P_{[a,b]}=\{([a,b],t)\}$ of $[a,b]$ and since $\lim\limits_{k \to \infty}f_k(x)=f(x)$ then there exists $m_2 > m_1 $ such that $S(f_{m_2}, P_{[a,b]}) -S(f,P) \in U$ then $S(f, P_{[a,b]})-\mathcal{A} \in U.$ Therefore $f$ is $ap-$Henstock-Kurzweil integrable on $[a,b],$ and $$\lim\limits_{k \to \infty}(ap)\int_{a}^{b}f_k  = (ap)\int_{a}^{b}f .$$
\end{proof}
\begin{definition}
\begin{enumerate}
\item[(1)] Let $\{f_k\},~\{{F}_k\}$ be sequences of functions defined on $[a,b]$ and let $E \subset [a,b]$ be a $\mu-$measurable set. Then $\{f_k\}$ is said to be uniformly $\mu_{AP}$-Henstock-Kurzweil integrable on $[a,b]$ if for every $\varepsilon>0$ there exists a (ANF) $S$ on $[a,b]$ such that 
\begin{eqnarray*}
S(f_k - F_k, P)-\mathcal{A} \in U
\end{eqnarray*}
for all $k$ whenever $P$ is $S-$tagged partitions where ${F}_k$ is the primitive of $f_k$ for each $k$ and $U$ is a $\theta-$nbd.  We denote it as $\lambda-UAP([a,b]).$
\item[(2)] $\{F_k\}$ is said to satisfy the uniformly approximate strong Lusin condition on $E$ (i.e., $F_k \in \lambda-ASL(E)$) if for every $E_1 \subset E$ with $\lambda(E_1)=0$ and for every $\varepsilon>0$ there exists an (ANF) $S$ on $E$ such that 
\begin{eqnarray*}
S(|F_k|,P)-\mathcal{A} \in U
\end{eqnarray*}
whenever $P$ belongs in a $S-$fine tagged partition $P_1$ of $E_1$ and $U$ is a $\theta-$nbd.
\end{enumerate}
\end{definition}
Now we discuss about $\mu_{AP}$-Henstock-Kurzweil equi-integrability and uniformly strong Lusin condition (in short ASL).
\begin{theorem}
Let $\{f_k\}$ be a sequence of functions $f_k: [a,b] \to X$ and let $f: [a,b] \to X$ be any function. If the followings are holds:
\begin{enumerate}
\item $\{f_k\} \to f(x)$ a.e. on $[a,b],$
\item $\{f_k\}$ is $ap$-Henstock-Kurzweil equi-integrable, 
\end{enumerate}
 this implies the followings are equivalent:
\begin{enumerate}
\item[(a)] $\{f_k\}$ is pointwise bounded,
\item[(b)] $\{F_k\}$ is $ASL.$
\end{enumerate}
\end{theorem}
\begin{proof}
The proof follows from the same technique as used in the \cite[Lemma 2.2]{Mema}.
\end{proof}
\begin{lemma}
Let $\{f_k\}$ be a sequence of measurable functions defined on $[a,b]$ satisfying the following conditions
\begin{enumerate}
\item $f_k(x) \to f(x)$ a.e. on $[a,b]$ as $k \to \infty.$
\item $\{F_k\} \in \lambda-ASL([a,b]),$ where $ F_k$ is the primitive of $f_k.$
\item $\{f_k\} \in \lambda-UAP([a,b]),$
\end{enumerate}
then $f \in AP([a,b], X)$ and $(ap)\int_{[a,b]}f = \lim\limits_{k \to \infty}(ap)\int_{[a,b]}f_k .$
\end{lemma}
\begin{proof}
The proof is similar as the \cite[Theorem 3.1]{Jae1}.
\end{proof}
\begin{theorem}
Let $\{f_k\}$ be a sequence of measurable functions on $[a,b]$ with the followings:
\begin{enumerate}
\item $\{f_k\}$ is pointwise bounded on $[a,b].$
\item $\{f_k\} \in \lambda-UAP([a,b]).$
\end{enumerate}
then $\{F_k\} \in \lambda-ASL([a,b]),$ where $F_k$ is the primitive of $f_k.$
\end{theorem}
\begin{proof}
Let $Y \subset [a,b]$ of $\lambda(Y)=0.$ Let $\varepsilon>0.$ For each $i,$ consider the set $Y_i=\{x \in Y:~i-1 \leq \sup\limits_{k}\lambda(f_k(x))<i\}.$ Choose an open set $O_i$ such that $Y_i \subset O_i$ and $\lambda(O_i)< \frac{\varepsilon_i}{i}.$ As $\{f_k\} \in \lambda-UAP([a,b])$ then there exists an ANF $S^\prime$ on $[a,b]$ and a $\theta-$nbd $U$  such that 
\begin{eqnarray*}
S(f_k-F_k,P)-(ap)\int_{[a,b]}f \in U 
\end{eqnarray*}
for all $k$ whenever $P$ is a $S-$fine partition of $[a,b].$ Let $\delta(x)>0$ on $Y_i$ so that $(x-\delta(x), x+ \delta(x)) \subset O_i$ when $ x \in Y_i.$ Let $S(x)$ be  defined on $[a,b]$ as 
\[
S(x)  = \left\{ {\begin{array}{*{20}c}
   S^\prime(x) \cap (x-\delta(x), x +\delta(x)~if~x \in Y_i,~i=1,2,..,  \\
   S^\prime(x)~if x \in [a,b] \setminus \bigcup Y_i  \\
 \end{array} } \right.
\] 
Let $P_i =\{([a,b],x) \in P:~ x \in Y_i\}$ then $P= \bigcup_{i}P_i$ where $P_i $ is in a $S-$fine partition as well as $S^\prime-$fine partitions. By using Saks-Henstock Lemma, we have
\begin{align*}
S(|F_k|,P)-\mathcal{A} \in U.
\end{align*}
So, $\{F_k\} \in \lambda-ASL([a,b]).$
\end{proof}
\begin{corollary}\label{cor11}
Let $\{f_k\}$ be a sequence of measurable functions defined on $[a,b]$ with the following conditions:
\begin{enumerate}
\item $f_k(x) \to f(x)$ a.e. on $[a,b]$ as $n \to \infty.$
\item $\{f_k\}$ is pointwise bounded on $[a,b].$
\item $\{f_k\} \in \mu-UAP(Q).$
\end{enumerate}
then $f \in AP([a,b], X))$ and $(ap)\int_{[a,b]}f  = \lim\limits_{ n \to \infty}(ap)\int_{[a,b]}f_k .$
\end{corollary}
From these above results we can find the following Theorem as:
\begin{theorem}
Let $\{f_k\}$ be a sequence of $\mu-$measurable  functions on $[a,b]$ with the following
\begin{enumerate}
\item $f_k(x) \to f(x)$ a.e. on $[a,b]$ as $k \to \infty.$
\item $\{f_k\}$ is uniformly bounded on $[a,b].$
\end{enumerate}
then $f \in AP([a,b],X))$ and $(ap)\int_{[a,b]}f  = \lim\limits_{ n \to \infty}(ap)\int_{[a,b]}f_k .$
\end{theorem}
\begin{proof}
As $|f_k(x)|\leq L$ for all $k$ and $x \in [a,b]$ with a positive constant $L.$ Since $f_k(x) \to f(x)$ a.e. on $[a,b]$ as $k \to \infty.$ This implies  $f_k$ and $f$ are  measurable and bounded a.e. on $[a,b]$,  hence $ap-$Henstock-Kurzweil integrable on $[a,b].$ Now from $(2),$ $\{f_k\}$ is uniformly bounded on $[a,b].$ Using the Corollary \ref{cor11}, we get the complete proof.

\end{proof}


\begin{thebibliography}{A}
\bibitem{B} R. G. Bartle, {\it A modern theory of integration}, Graduate Studies in Mathematics, Vol. 32,
American Mathematical Society, Providence, 2001.
\bibitem{Burkill}J. C. Burkill, The approximately continuous Perron integral, Mathematische
Zeitschrift 34 (1932), 270-–278.

\bibitem{HS} R. Henstock, {\it The General Theory of Integration} (Clarendon
Press, Oxford, 1991)
\bibitem{HI} T. H. Hildebrandt, On bounded functional operations. Trans. Amer. Math. Soc. 36(1934), 868--875.
\bibitem{J} F. Jones, Lebesgue Integration on Euclidean Space, Revised Edn.
(Jones and Bartlett Publishers, Boston, 2001)
\bibitem{Base}  J. Kurzweil, J. Jarnik, A convergence theorem 
for Henstock-Kurzweil integral and its relation to topology, Bulletin de la Classe des Sciences,  8(7-12)(1997)  217--230. 
\bibitem{VM}V. Marraffa, Riemann type integrals for functions taking values in a locally convex space, Czechoslovak Mathematical Journal, 56 (131) (2006), 475--489.
\bibitem{Mema} E. Mema, Equiintegrability and Controlled Convergence for the Henstock-Kurzweil Integral, International Mathematics Forum,  8(19)(2013), 913--919. 
\bibitem{KW} J. Kurzweil, Nichtabsolut konvergente Integrale. Teubner-Texte
z\"ur Mathematik, Band 26 (Teubner Verlagsgesellschaft, Leipzig,
1980).
\bibitem{O} K. M. Ostaszewski, Henstock integration in the plane, Memoirs of the American
Mathematical Society 353, American Mathematical Society, Providence, 1986.
\bibitem{Jae} J. M. Park, J.J. Oh, Chun-Gil Park, D. H. Lee, The AP-Denjoy and AP-Henstock integrals, Czechoslovak Mathematics Journal, 57(2)(2007), 689--696.
\bibitem{Jae1}J. M. Park, D. H. Lee, J. H. Yoon,  B. M. Kim, The convergence theorem for the $AP-$integral, Journal of Chungcheong Mathematical  Society,  12(1)(1999), 113--118.
\bibitem{RA}  R. A. Gordon, \textit{The Integrals of Lebesgue, Denjoy, Perron and Henstock}, Graduate Studies in Mathematics Vol 4, American Mathematical Society, (1994).
\bibitem{Val1} V. A. Skvortsov, P. Sworowski, The AP-Denjoy and AP-Henstock integrals revisited, Czechoslovak Mathematics Journal, 62(3)(2012), 581--591.
\bibitem{Val} V. A. Skvortsov, T. Sworowska, P. Sworowski, On approximately continuous integrals (A survey), Traditional and present-day topics in real analysis, Wydawnictwo Uniwersytetu Lodzkiego/Lodz University Press (Poland) (2013), https://DOI: 10.18778/7525-971-1.15.
\bibitem{Kwan} K. Shin, J.H. Yoon, Properties for $AP-$Henstock integral, Reprint, 2021.
\bibitem{JU} J. H. Yoon, The AP-Henstock integral of vector-valued functions, Journal of Chuncheong Mathematical Society, 30(1)(2017), 151--157.


 \end{thebibliography}
\end{document}